\documentclass[12pt,leqno]{article}
\usepackage{amsmath,amssymb,amsthm,amscd}
\let\ftilde\tilde 
\renewcommand{\tilde}[1]{\mathbf{\ftilde{\mathnormal{#1}}}}

\newcommand{\C}{\mathbb{C}}
\newcommand{\Z}{\mathbb{Z}}

\newcommand{\unit}{\mathbf{1}}

\newcommand{\tr}{\operatorname{tr}}

\newcommand{\frakg}{\mathfrak{g}}
\newcommand{\frakh}{\mathfrak{h}}

\newcommand{\wt}[1]{\operatorname{wt}\,(#1)}
\newcommand{\ch}{\operatorname{ch}}
\newcommand{\NO}{\,{\raise0.25em\hbox{$\mathop{\hphantom{\cdot}}%
\limits^{_{\circ}}_{^{\circ}}$}}\,}

\def \tr{{\rm tr}}
\def \Res{{\rm Res}}

\def \<{\langle}
\def \>{\rangle}

\def \l{\lambda }

\renewcommand{\theequation}{\thesection.\arabic{equation}}
\catcode`\@=11
\renewcommand{\theequation}{%
\thesection.\arabic{equation}}
\@addtoreset{equation}{section}
\makeatother
\catcode`\@=\active

\theoremstyle{plain}
 \newtheorem{thmm}{Theorem}
 \newtheorem{theorem}{Theorem}[subsection]
 \newtheorem{corollary}[theorem]{Corollary}
 \newtheorem{lemma}[theorem]{Lemma}
 \newtheorem{proposition}[theorem]{Proposition}
 
\theoremstyle{remark}

\begin{document}
\begin{Large}
\begin{center}
\textbf{Quasi-modular forms and trace functions associated to
free boson and lattice vertex operator algebras}
\end{center}
\end{Large}

\vskip 5ex

\begin{center} Chongying Dong\footnote{Supported by NSF grant 
DMS-9700923 and a research grant from the Committee on Research, 
UC Santa Cruz.}, Geoffrey Mason\footnote{Supported by NSF grant DMS-9700909 and a research grant from the Committee on Research, UC Santa Cruz.}\\ Department of
Mathematics, University of California  Santa Cruz, CA 95064, U.S.A\\

\vskip 2ex Kiyokazu Nagatomo\footnote{Supported in part by  Grant-in-Aid
for Scientific Research, the Ministry of Education, Science and
Culture.}\\ Department of Mathematics, Graduate School of Science, Osaka
University\\ Osaka, Toyonaka 560-0043, Japan
\end{center}

\vskip 10ex


\begin{small}

\noindent
\textbf{Abstract}: We study graded traces of vectors in free bosonic
vertex operator algebras and lattice vertex operator algebras. 
We show in particular that trace functions in these
two theories always have the shape $f(q)/\eta(q)^d$ where $f(q)$ is
quasi-modular in the case of $d$ free bosons, and modular (i.e., a sum of 
holomorphic modular forms of various weights) in the case of theories
based on a lattice $L$ of rank $d.$ We also show how spherical harmonic
polynomials with respect to $L$ are related to primary fields in lattice theories.
\end{small}

\section{Introduction}

The purpose of this paper is to study the space of 1-point 
correlation functions, or trace functions, which arise from certain 
vertex operator algebras. 

Let us suppose that $V$ is a vertex operator algebra \cite{FLM} (also see
\cite{B}, \cite{MN}) with standard
$L(0)$-grading 
\[
V =\bigoplus_{n\in\Z} V_n.
\]
The most basic 
trace function is the formal graded character
\[
\ch_q V =\tr_V q^{L(0) - c/24} =  q^{-c/24}\sum_{n\in\Z}(\dim V_n)q^n
\]
where $c$ is the central charge of $V.$

      For many well-known VOAs, the graded character has certain 
modular-invariance properties. For example, if  $V(\frakg,l)$ is the 
generalized Verma module of level $l$ associated to a (complex) Lie 
algebra $\frakg$ of dimension $d$ (cf.  [FZ], [L]) and 
trivial $\frakg$-module, then $V(\frakg,l)$ is a 
vertex (operator) algebra satisfying
\[
\ch_q V(\frakg,l) = \eta(q)^{-d}
\]
where $\eta(q)$ is the Dedekind eta-function. Similarly, if $L$ is a 
positive-definite, even lattice of rank $d$ and  $\theta_L(q)$
is the corresponding theta-function, 
and if $V_L$ is the associated VOA (cf. [B], [FLM]), then
\[
\ch_q V_L =\theta_L(q)/ \eta(q)^d.
\]
So $\ch_q V_L$  is a modular function (i.e., of weight zero) on 
some congruence subgroup of the modular group $SL(2,\Z),$ while $\ch_q 
V(\frakg,l)$ is a modular form of weight
$-d/2$ on the full modular group\footnote{We will have no need to concern 
ourselves with the fact that $d$ may be odd, so that the weight may be a half 
integer.}.  On the other hand, the VOA $M(c,0)$ 
associated to the Virasoro algebra of central charge $c$ (cf. [FZ], [L]) 
satisfies
\[
\ch_q M(c,0) = (1-q)q^{-c/24+1}/\eta(q)
\]
and hence is not modular. Here, of course, we are interpreting $q$ in 
the usual way to be
equal to $e^{2 \pi i \tau}$ with $\tau$ in the complex upper-half plane $H.$

      In the fundamental paper of Zhu [Zh] general correlation 
functions, which generalize 
the graded character, were studied. Recall (loc cit) that if $v$ lies 
in $V_k$  with vertex operator $Y(v,z) = \sum_{ n \in\Z}v(n) z^{-n-1},$ then 
the so-called zero mode $o(v) = v(k-1)$ is a linear operator on V which 
leaves invariant each homogeneous space $V_n,$ so that one can form the 
expression
\[
Z(v, q) = \tr_V o(v)q^{L(0)-c/24}=q^{-c/24} \sum_{n\in\Z}(\tr_{V_n} o(v))q^n.
\]

     One calls the linear extension of $Z(v,q)$ to all of V the 
(1-point) correlation function determined by $V.$ It is the purpose of 
this paper to understand the nature of this function   in the case of 
the lattice VOAs $V_L$ and the Heisenberg VOA $M(1),$ which is the space 
$V(\frakg,l)$ in the case that $\frakg$
 is an abelian Lie algebra equipped with a 
non-degenerate symmetric bilinear form $(,)$ and $l=1.$ $M(1)$ is often 
referred to as the VOA of $d$ free bosons if $\dim \frakg = d.$ In particular, 
we ask when $Z(v,q)$ is modular in a suitable sense and, if not, how 
does it deviate from being modular?

     In order to explain our results and methods, we need to introduce
so-called quasi-modular forms [KZ], which play an interesting role in
the present paper. For a positive integer $N$ and for a Dirichlet
character $\epsilon$ modulo $N$, let $M(N,\epsilon)$ denote the ring of
modular forms on the congruence subgroup $\Gamma_0(N)$ which transform
according to the character $\epsilon$ and which are holomorphic in $H.$ 
Let
\begin{equation}\label{1.0}
Q(N,\epsilon) = M(N,\epsilon) [E_2]
\end{equation}
be the space obtained by adjoining to $M(N,\epsilon)$ the Eisenstein 
series $E_2(q).$ For example, if $N = 1,$ so that also $\epsilon = 1,$ then 
we have
\begin{equation}\label{1.1}
Q =  Q(1,1) = C[E_2, E_4,E_6]
\end{equation}
the full ring of quasi-modular forms on $SL(2,\Z).$ We can now state our 
first result:
\begin{thmm}\label{thmm1}
Let $V$ be the VOA of $d$ free bosons, with $Q$ as in (\ref{1.1}).

(a) If $v$ is in $V,$ then $Z(v,q)$ converges to a holomorphic function
$f(v,q)/\eta(q)^d$ in the upper half plane 
for some  $f(v,q)\in Q.$

(b) Every $f(q)$ in $Q$ may be realized as $f(v,q)$ for some $v$ in $V.$
\end{thmm}

     We can restate Theorem \ref{thmm1} in the following form: there is a 
linear surjection
\begin{equation}\label{1.2}
\begin{array}{cccl}
t:  &V& \to &Q \\
   &v& \mapsto &(\ch_qV)^{-1}Z(v, q). 
\end{array}
\end{equation}
In fact we will construct a section of the map t, more precisely we 
explicitly describe a subspace $W$ of $V$ such that the restriction of $t$ 
to $W$ is an isomorphism onto $Q.$

\begin{thmm}\label{thmm2} Let $L$ be a positive-definite even lattice of 
rank $d$ and level $N,$ so that the theta function $\theta_L(q)$ lies in 
$M(N,\epsilon)$ for suitable $\epsilon.$ Then for every element $v$ in $V_L,$
 we have
\begin{equation}\label{1.3}
Z(v,q) =f(v,q)/\eta(q)^d
\end{equation}
for some $f(v,q)$ in $M(N,\epsilon).$ In particular, each $Z(v,q)$ is a sum 
of modular forms (of varying weights).
\end{thmm}

    The role of quasi-modular forms in the proof of Theorem
\ref{thmm2} is more hidden. To explain, we recall here (and discuss in
greater detail in Section 4.1) that in [Zh], Zhu shows how to identify
vectors in $V$ such that the corresponding trace function $Z(v,q)$
essentially behaves like a modular form of weight $k$ for some $k.$ 
More precisely, if we think of $Z(v,q)$ as a function of
$\tau,$ then under the action of the modular group
\[
\gamma =\left(\begin{array}{cc} a & b\\ c & d\end{array}\right):  
Z(v, \tau)  \mapsto (c \tau +d)^{-k}Z(v, \gamma\tau),\ \gamma\in SL(2,\Z),
\]
the $\gamma$-transforms of $Z(v,\tau)$ spans a finite-dimensional $SL(2,\Z)$-module of holomorphic functions on $H.$ Moreover
each such $\gamma$ transform has a suitable q-expansion. Extending the 
language of [KM], we call such functions {\em generalized modular forms}
of weight $k.$ A 
modular form has this property, however it is shown in [KM] that 
there are generalized modular forms which are not modular. 
Thus one cannot deduce from 
Zhu's results alone that the trace functions $Z(v,q)$ are modular, or 
even sums of modular forms. However, we will show below that the 
trace functions occurring in Theorem \ref{thmm2} also have the property that
$\eta(q)^dZ(v,q)$ lies in $Q(N,\epsilon)$ 
so that in this regard they behave similarly to the trace functions 
in the free boson case. It is easy to see that a generalized modular 
form that is also a quasi-modular form is necessarily an ordinary 
modular form, and Theorem \ref{thmm2} then follows.

    Our third main theorem has a slightly different flavor. In the 
paper [DM1] the space of trace functions $Z(v,q)$ was determined in the 
case of the Moonshine module $V^{\natural}$ [FLM],  a prominent role being 
played by the primary fields i.e., the (homogeneous) vectors of $V^{\natural}$ 
which are highest weight vectors for the Virasoro algebra. Based on 
this work and the calculations in [HL], we conjecture that if $V$ is a 
holomorphic VOA then each cusp form on $SL(2,\Z)$ (possibly with character)
can be realized by a trace function $Z(v,q)$ in which $v$ is a primary 
field.  The conjecture seems to be 
non-trivial for any  holomorphic VOA. We will prove it for the 
lattice VOA $V_{E_8}$ based on the $E_8$ root lattice. Our approach 
uses the theory of spherical harmonics, and through this mechanism 
one sees that the conjecture may be viewed as a conformal field 
theoretic analog of the following problem in number theory: given a 
positive-definite self-dual even lattice $L,$ describe the space of 
modular forms $\theta_L(P,\tau)$ obtained by modifying $\theta_L(\tau)$ by 
$P,$ where $P$ ranges over the homogeneous spherical harmonic functions 
with respect to $L.$ (Replacing
``holomorphic" by ``rational" in our conjecture corresponds to 
eliminating the self-duality of $L.$) Waldspurger showed [W] that all 
cusp forms of level one can be obtained as $\theta_{E_8}(P,\tau)$ for 
suitable $P,$ and this leads to our result about $V_{E_8}$ because of the 
following.

\begin{thmm}\label{thmm3} Let $P$ be a homogeneous spherical harmonic of 
degree $k$ with respect to the lattice $L$ of rank $d.$ Then there is a 
primary field $v_P$ in the lattice VOA $V_L$ with the property that
\[
Z(v_P, q) = \theta_L(P,q)/ \eta(q)^d.
\]
\end{thmm}

    There is a circle of ideas relating our results: in [EZ], Eichler 
and Zagier gave an   approach to the result of Waldspurger using 
ideas from the theory of Jacobi forms. Jacobi forms and Jacobi-like 
forms [Za] are closely related to quasi-modular forms, and they also 
play a role in the results of [DM2] which identify the quasi-modular 
forms underlying the proof of theorem 2 and relate them to theta 
functions modified by a spherical harmonic. And the occurrence of 
Jacobi forms in string theory has been observed frequently in the 
past few years.

\section{Vertex operator algebras and graded traces}

In this section we briefly review from [Zh] the ``bracket''  vertex operator 
algebra $(V,Y[\ ], \unit, \omega-c/24)$ constructed 
from $(V,Y,\unit, \omega).$ 
We also list several formulas involving  the trace functions
from [Zh] and present an easy corollary which is used
frequently in the later sections. 

\renewcommand{\theequation}{\thesubsection.\arabic{equation}}
\catcode`\@=11
\renewcommand{\theequation}{%
\thesubsection.\arabic{equation}}
\@addtoreset{equation}{subsection}
\makeatother
\catcode`\@=\active

\subsection{Vertex operator algebras of genus one}
\label{section:2.1}

Let $V=(V,Y,\unit,\omega)$ be a vertex operator algebra. Then $V$ may be
regarded as a vertex operator algebra on the sphere. In order to study
modular invariance in the theory of vertex operator algebras,
a new vertex operator algebra on the torus was introduced in \cite{Zh}. 
The new vertex operator algebra is $(V,Y[\ ],\unit, \omega-c/24)$ where 
$c$ is the
central charge of $V.$ The new vertex 
operator associated to  a homogeneous element $a$ is given by
\[ Y[a,z] = \sum_{n\in\Z}a[n]z^{-n-1} = Y(a,e^{z} -1)e^{z\wt{a}}
\] while a Virasoro element is $\tilde{\omega} = \omega-c/24$.
Thus 
\[
a[m] = \Res_z\left(Y(a,z)(\ln{(1+z)})^m(1+z)^{\wt{a}-1}\right)
\]
and \[
a[m] = \sum_{i=m}^\infty c(\wt{a},i,m)a(i)
\]
for some scalars $c(\wt{a},i,m)$ such that $c(\wt{a},m,m)=1.$ 
In particular,
\[
a[0]=\sum_{i\geq 0}\binom{\wt{a}-1}{i}a(i).
\]
We also write
\[
L[z] = Y[\omega,z] = \sum_{n\in\Z} L[n]z^{-n-2}.
\]
Then the $L[n]$ again generate a copy of the Virasoro algebra with
the same central charge $c.$ Now $V$ is graded by
the $L[0]$-eigenvalues, that is
\[
V=\bigoplus_{n\in\Z}V_{[n]}
\]
where $V_{[n]}=\{v\in V|L[0]v=nv\}.$ We also write $\operatorname{wt}[a]=n$
if $a\in V_{[n]}.$ It should be pointed out that for any
$n\in\Z$ we have 
\[
\sum_{m\leq n}V_n=\sum_{m\leq n}V_{[n]}.
\]

We also recall the notion of $V$-module briefly. A $V$-module $M=(M,Y^M)$
is a $\C$-graded vector space 
\[
M=\oplus_{\lambda\in \C}M_{\lambda}
\]
such that each $M_{\l}$ is finite-dimensional and $M_{\l+n}$ is zero
if n is a small enough integer. Furthermore $M_{\l}$ is the eigenspace
of $L(0)$ with eigenvalue $\l$ where $L(0)$ is the component
operator of $Y^M(\omega,z)=\sum_{n\in\Z}L(n)z^{-n-2}.$ As in the
case of vertex operator algebras, the component operator
$o(a)=a^M(\wt{a}-1)$ preserves each homogeneous space $M_{\l}$ if
$a\in V$ is homogeneous and $Y^M(a,z)=\sum_{n\in\Z}a^M(n)z^{-n-1}.$
Again using the $\C$-linearity, we extend the operator $o(a)$ to
all $a\in V.$

\subsection{Graded traces}
\label{section:2.2}

Next we  state some results from \cite{Zh}. We use the Eisenstein series
$E_{2k}(\tau)$ normalized as in \cite{DLM} equation (4.28). Thus:
\[
E_{2k}(\tau) =\frac{-B_{2k}}{2k!}+\frac{2}{(2k-1)!}
\sum_{n=1}^\infty \sigma_{2k-1}(n)q^n
\]
where $\sigma_k(n)$ is the sum of the $k$-powers of the divisors of
$n$ and $B_{2k}$ a Bernoulli number.

Let $M=(M,Y^M)$ be a $V$-module. For any $a\in V$ we define
a formal power series in $q:$ 
\[
Z_M(a,q)=\tr_Mo(a)q^{L(0)-c/24}=q^{-c/24}\sum_{\l\in\C}(\tr_{M_{\l}}o(a))q^{\l}.
\]
Then $Z_V(a,q)$ is exactly $Z(a,q)$ defined before.

\begin{proposition}[\cite{Zh}, Proposition 4.3.5, 4.3.6]
\label{proposition:2.2.1} 
Let $M$ be a $V$-module. 
Then the following identities hold as formal power series for any $a,b\in V.$
\begin{align}
&Z_M(a[0]b,q)= 0,\label{eqn:2.2.1}\\
&\tr_M\,o(a)o(b)q^{L(0)-c/24} = Z_M(a[-1]b,q)-
\sum_{k=1}^\infty{E}_{2k}(q)Z_M(a[2k-1]b,q),\label{eqn:3.2.1}\\
&Z_M(a[-2]b,q)=-\sum_{k=2}^\infty(2k-1){E}_{2k}(q)
Z_M(a[2k-2]b,q).\label{eqn:2.2.3}
\end{align}
\end{proposition}

The following corollary plays an important role in our arguments.

\begin{corollary}
\label{corollary:2.2.2} Let notation  be as before. Then 
for each positive integer $r$ we have
\begin{align}\label{00}
Z_M(a[-r]b, q)&=\delta_{1,r}\tr_Mo(a)o(b)q^{L(0)-c/24}\\
&\ \ \ \ +(-1)^{r+1}\sum_{k>r/2}^{\infty}h(k,r){E}_{2k}(q)Z_M(a[2k-r]b,q)\nonumber
\end{align}
where $h(k,r) = \binom{2k-1}{r-1}.$
\end{corollary}
\begin{proof}
We  prove the result by induction on $r$. The cases $r=1,2$ are nothing
but (\ref{eqn:3.2.1}) and (\ref{eqn:2.2.3}) respectively. Suppose the statement is
true for $r\geq 2$. Then replacing $a$ by $L[-1]a$ and noting
$(L[-1]a)[n] = -na[n-1]$, we see that
\[
rZ_M(a[-r-1]b,q)=(-1)^{r+1}\sum_{k> r/2}
(r-2k)h(k,r){E}_{2k}(q)Z_M(a[2k-r-1]b,q).
\]
Then by (\ref{eqn:2.2.1}), we conclude that
\[
Z_M(a[-r-1]b,q)=(-1)^{r+2}\sum_{k>(r+1)/2}^\infty
h(k,r+1){E}_{2k}(q)Z_M(a[2k-r-1]b,q).
\]
\end{proof}

\section{Quasi-modularity of trace functions for free boson VOAs}
\label{section:3}

In this section,  after reviewing from [FLM] 
the free boson vertex operator 
algebra $M(1),$ we establish Theorem \ref{thmm1}.

\subsection{Free boson VOAs}
\label{section:3.1}

Let $\frakh$ be a $d$-dimensional vector space with a non-degenerate
symmetric bilinear form $(\,,\,)$ and $\hat{\frakh}$ be the corresponding
affinization viewing $\frakh$ as an abelian Lie algebra:
$\hat{\frakh} =\frakh\otimes \C[t,t^{-1}]\oplus\C K$ with commutator
relations
\begin{align*} &[h\otimes t^m,h'\otimes t^n] =
(h,h')\delta_{m+n,0}K,\quad (h,h'\in\frakh, m,n\in \Z),\\
&[K,\frakh\otimes
\C[t,t^{-1}]] = 0.
\end{align*} Consider the induced module
\[
M(1) = U(\hat{\frakh})\otimes_{\frakh\otimes \C[t]\oplus \C K}\C
\] where $\frakh\otimes \C[t]$ acts trivially on $\C$ and $K$ acts as $1$.
We denote by $h(n)$ the action of $h\otimes t^n$ on $M(1)$. The space
$M(1)$ is linearly isomorphic to the symmetric algebra
$S(\frakh\otimes t^{-1}\C[t^{-1}])$.  Thus setting
$\unit = 1\otimes 1$, any element in $M(1)$ is a linear
combination of elements of type
\[ v = a_1(-n_1)\cdots a_k(-n_k)\unit,\,(a_1,\dots,a_k\in
\frakh, n_1,\dots,n_k\in \Z_{+}).
\] For such $v$, we define
\begin{equation}\label{d1}
Y(v,z) =
\NO\partial^{(n_1-1)}a_1(z)\cdots\partial^{(n_k-1)}a_k(z)\NO,
\quad \partial^{(n)} = \frac{1}{n!}\left(\frac{d}{dz}\right)^n
\end{equation}
 where
\begin{equation}\label{d2}
a(z)=\sum_{n\in\Z}a(n)z^{-n-1}
\end{equation}
and 
$\NO\,\NO$ indicates the normal ordering procedure.

Now let $\{h_i\}_{i=1}^d$ be an orthonormal basis of $\frakh$ and set
$\omega = \frac{1}{2}\sum_{i=1}^d h_i(-1)^2\unit$. Then $(M(1), Y,
\unit,\omega)$ is a vertex operator algebra with a vacuum $\unit$ and 
Virasoro element $\omega$ (see \cite{FLM}). In particular,
\[
M(1)=\bigoplus_{n\geq 0}M(1)_n
\]
where $ M(1)_n=\<a_1(-n_1)\cdots a_k(-n_k)\unit|a_1,\dots,a_k\in
\frakh, n_1,\dots,n_k\in \Z_{+}, \sum n_i=n\>.$
We will identify $M(1)_1$ with $\frakh$ in an obvious way.

Note that $o(a)=0$ for $a\in \frakh,$ so that (\ref{00}) simplifies
in the free boson case. Moreover, we see from Subsection 2.1 that 
$a[0]=a(0)=0.$  Thus we have the following commutator relation
by noting that $a[1]b=a(1)b=(a,b)$ for $a,b\in  \frakh:$
\begin{equation}\label{ae}
[a[m], b[n]]=m\delta_{m+n,0}(a,b).
\end{equation}
It is easy to see from this that $M(1)_{[n]}$ is spanned by vectors
 $a_1[-n_1]\cdots a_k[-n_k]\unit$ for $a_1,\dots,a_k\in
\frakh, n_1,\dots,n_k\in \Z_{+}$ such that $\sum n_i=n.$

\subsection{Proof of Theorem 1}
\label{section:3.3}

We first prove part (a) of Theorem \ref{thmm1}.  It is enough 
to prove it for an element of the spanning set defined in Section 3.1.

Let $v=a_1[-r_1]\cdots a_k[-r_k]\unit$ and
define the length of $v$ to be $k,$ denoting this by $l(v)=k.$ 
 We prove the statement by induction on $k.$ If $k=0$ then $v=\unit,$
$Z(v,q)=1/\eta(q)^d$ and there is nothing to prove.
 Now we assume that the statement holds for
$v$ with $l(v)\leq k,$ and consider $v=a[-r]a_1[-r_1]\cdots
a_k[-r_k]\unit=a[-r]b$ where $b=a_1[-r_1]\cdots a_k[-r_k]\unit.$ 
Since $o(a)=0$, we see from equation (\ref{00}) that
\begin{equation}\label{01}
Z(a[-r]b, q)=\sum_{m>r/2}^{\infty}h(m,r){E}_{2m}(q) Z(a[2m-r]b,q).
\end{equation}
Since $2m-r>0,$ $a[2m-r]b$ is a linear combination of homogeneous
vectors whose lengths are no greater than $k.$ By induction,
each $Z(a[2m-r]b,q)$ converges to a holomorphic function which is a
quotient of a quasi-modular form by $\eta(q)^d.$ Note that
${E}_{2m}(q)$ is in $Q.$   It is immediate
that $Z(v,q)$ converges to a holomorphic function of the same kind.

In order to prove part (b) of Theorem \ref{thmm1} we need several lemmas.
From now on, we fix $a\in \frakh$ such that $(a,a)=1.$ 
\begin{lemma}\label{l1} If $n>0,$  $r\geq 0$ are integers,  we have
\[
Z(a[-n]^{2r}\unit,q)=(-1)^{(n+1)r}n^r(2r-1)!!h(n,n)^rE_{2n}(q)^r/\eta(q)^d
\]
where $(2r-1)!!=1\cdot 3\cdot 5\cdots (2r-3)\cdot (2r-1).$ 
\end{lemma}
\begin{proof} Recall Corollary \ref{corollary:2.2.2} and the fact that
$a[s]a[-t]=s\delta_{s,t}$ for positive integers $s,t$ (cf. (\ref{ae})).
We have 
\begin{align*}
Z(a[-n]^{2r}\unit, q)&=(-1)^{n+1}h(n,n)E_{2n}(q)Z(a[2n-n]a[-n]^{2r-1}\unit,q)\\
&=(-1)^{n+1}n(2r-1)h(n,n)E_{2n}(q)Z(a[-n]^{2r-2}\unit,q).
\end{align*}
Now the result follows by induction, noting that $Z(\unit,q)=1/\eta(q)^d.$
\end{proof}

\begin{lemma}\label{l2} Let $r,s\geq 0$ be integers. Then 
\begin{equation}\label{03}
Z(a[-1]^{2r}a[-2]^{2s}\unit,q) =(-6)^s(2r-1)!!(2s-1)!!E_2(q)^r E_4(q)^s/\eta(q)^d
\end{equation}
and
\begin{equation}\label{04}
Z(a[-2]^{2s}a[-3]^{2t}\unit,q) =(-6)^s(2s-1)!!(2t-1)!!(30)^tE_4(q)^s E_6(q)^t/\eta(q)^d.
\end{equation}
\end{lemma}
\begin{proof}
Since the proofs of (\ref{03}) and (\ref{04}) are almost identical we only
prove (\ref{03}). The special case $r=0$ follows from 
Lemma \ref{l1}. So we can proceed by induction on $r$ and assume that $r\ne 0.$ Then by Corollary
\ref{corollary:2.2.2} and Lemma \ref{l1} we have
\begin{align*}
Z(a[-1]^{2r}a[-2]^{2s}\unit,q)&=(2r-1)h(1,1)E_2(q)Z(a[-1]^{2r-2}a[-2]^{2s}\unit,q)\\
&=(-6)^s(2r-1)!!(2s-1)!!E_2(q)^rE_4(q)^s/\eta(q)^d,
\end{align*}
as required.
\end{proof}

\begin{lemma}\label{l3} Let $r,s,t\geq 0$ be integers. Then  
\begin{equation}\label{02}
Z(a[-1]^{2r}a[-2]^{2s}a[-3]^{2t}\unit,q) =
c_{r,s,t}E_2(q)^r E_4(q)^s E_6(q)^t/\eta(q)^d+ lower\ terms
\end{equation}
where $c_{r,s,t}$ is a nonzero constant 
and ``lower terms" means a linear combination  of terms 
of the form $E_2(q)^{r'} E_4(q)^{s'}E_6(q)^{t'}/\eta(q)^d$
with $0\leq r'<r$ if $r\geq 1,$ and $0$ if $r=0.$ 
\end{lemma}

\begin{proof} If $r=0$ this is (\ref{04}), so we may assume that $r\geq 1.$
First we deal with the case $s=0.$ Following 
Corollary \ref{corollary:2.2.2} and Lemma \ref{l1} we have 
\begin{align*}
Z(a[-1]^{2r}a[-2]^{2s}a[-3]^{2t}\unit,q)&=(2r-1)E_2(q)Z(a[-1]^{2r-2}a[-2]^{2s}a[-3]^{2t}\unit,q)\\
&\ \ \ +6tE_4(q)Z(a[-1]^{2r-1}a[-2]^{2s}a[-3]^{2t-1}\unit,q).
\end{align*}
Continuing in this fashion, the result follows.
\end{proof}

Now we can prove part (b) of Theorem \ref{thmm1}. Since the space $Q$ of
quasi-modular forms on $SL(2,\Z)$ is an algebra generated by $E_2(q), E_4(q)$ and $E_6(q),$
it is enough to prove (b) for $f(q)=E_2(q)^rE_4(q)^sE_6(q)^t$ with
 non-negative integers $r,s,t.$ But this follows from Lemmas
\ref{l2} and \ref{l3}. Indeed, this shows that if $W$ is the subspace
of $M(1)$ spanned by $a[-1]^{2r}a[-2]^{2s}a[-3]^{2t}\unit$ for
$r,s,t\geq 0$ then the map $t$ of (\ref{1.2}) induces an isomorphism from
$W$ to $Q.$

\section{Modularity of trace functions for lattice vertex operator algebras}

We study  trace functions for lattice vertex operator
algebras $V_L$ ([B], [FLM]). We prove that
$Z(v,q)$ is quasi-modular for all $v\in V_L.$ This result
together with a result in [DM2] yields a proof of Theorem \ref{thmm3}.
It should be mentioned that in fact, we compute the
trace function $Z_M(v,q)$ corresponding to an irreducible module $M$. It is 
clear from our proof that each $Z_M(v,q)$ is also modular.

\subsection{Vertex operator algebras associated to lattices}

We now work in the setting of Chapter 8 of \cite{FLM}. In particular,
$L$ is a positive definite lattice, and $\frakh=L\otimes_{\Z}\C.$ Let
$\C[L]$ be the group algebra with a basis 
$\{e^{\alpha}|\alpha\in L\}.$ Then the vertex operator algebra
$V_L$ associated to $L$ is
\[
V_L=M(1)\otimes \C[L]
\]
as vector spaces. The vertex operator $Y(v,z)$ for $v\in M(1)$ is
defined as in the case of $M(1)$ (see formulas (\ref{d1}) and (\ref{d2})).
The operator $a(n)$ for $a\in\frakh$ and $n\ne 0$ acts on $V_L$ via its action
on $M(1).$ The operator $a(0)$ acts on $V_L$ by acting on $\C[L]$ in the following way:
\begin{equation}\label{4.1.1}
a(0)e^{\alpha}=(a,\alpha)e^{\alpha}
\end{equation}
for $\alpha\in L.$ We identify $M(1)$ with $M(1)\otimes e^0.$ Then 
the vacuum of $V_L$  is the vacuum  $\unit$ of $M(1)$ and the Virasoro element
is the same as before. For the purposes of this paper we do not
make explicit the expression for vertex operators $Y(v,z)$ for general $v$ 
except for those $v\in M(1).$ We remark that $V_L$ is a   
rational
vertex operator algebra in the sense of \cite{DLM}.

Let $L^{\circ}=\{\alpha\in L\otimes_{\Z}\C|(\alpha,L)\subset\Z\}$ be
the dual lattice of $L$ with coset decomposition $L^{\circ}=\cup_{i\in
L^{\circ}/L}(L+\l_i).$ Then each space $V_{L+\l_i}=M(1)\otimes
\C[L+\l_i]$ is an irreducible $V_L$-module \cite{FLM} where
$\C[L+\l_i]$ is the subspace of the group algebra $\C[L^{\circ}]$
corresponding the coset $L+\l_i.$ Moreover $\{V_{L+\lambda_i}|i\in
L^{\circ}/L\}$ constitute the complete set of non-isomorphic
irreducible $V_L$-modules \cite{D}. Since $V_L$ also satisfies the
$C_2$-condition (cf. \cite{Zh}, \cite{DLM}) we have the following
modular property for $V_L$ by \cite{Zh}:
\begin{proposition}\label{pa} Let $v\in (V_L)_{[n]}$ and $\gamma=\left(\begin{array}{cc}
a & b\\
c & d 
\end{array}\right)\in SL(2,\Z).$ Then $Z_i(v,q)$ 
converges to a holomorphic function in the upper half plane
and there exist scalars $c^{\gamma}_{ij}$
independent of $v$ and $\tau$  such that
\[
Z_i(v,\frac{a\tau+b}{c\tau+d})=(c\tau+d)^n\sum_{j\in L^{\circ}/L}c^{\gamma}_{i,j}Z_j(v,\tau)
\]
where $Z_i(v,\tau)=Z_i(v,q)=Z_{ V_{L+\l_i}}(v,q).$
\end{proposition}

\subsection{Graded traces in $V_L$}

In this subsection we determine $Z_i(v,q)$ for $i\in L^{\circ}/L$ and
$v\in V_L,$ and prove Theorem \ref{thmm2}. 

\begin{lemma}\label{l3.1} Let $v\in M(1)\otimes e^{\alpha}$ for nonzero
$\alpha\in L.$ Then $Z_i(v,q)=0.$
\end{lemma}

\begin{proof} Recall from \cite{FLM} the vertex operator $Y(v,z).$ It is
easy to see that $v(n)M(1)\otimes e^{\beta}\subset M(1)\otimes e^{\alpha+\beta}$ for $\beta\in L+\l_i$ for all $n\in\Z.$ So it is immediate that
$Z_i(v,q)=0.$
\end{proof}

It remains to determine $Z_i(v,q)$ for $v\in M(1).$ Note that as 
$M(1)$-module, $V_{L+\l_i}=\oplus_{\alpha\in L+\l_i}M(1)\otimes e^{\alpha}$
and each $M(1)\otimes e^{\alpha}$ is an irreducible $M(1)$-module.
For any $v\in M(1)$ and $\alpha\in L^{\circ}$ we set
\[
Z_{\alpha}(v,q)=Z_{M(1)\otimes e^{\alpha}}(v,q).
\]
Then 
\[
Z_i(v,q)=\sum_{\alpha\in L+\l_i}Z_{\alpha}(v,q).
\]
\begin{lemma}\label{l3.2} Let $a\in\frakh$ such that $(a,a)=1$ and $\alpha\in L^{\circ}.$ Then for any non-negative integer $r$, $Z_{\alpha}(a[-1]^r\unit,q)$
converges to holomorphic function in the upper half plane. Moreover
 there exist scalars $c_{r,r-2i}$ with $0\leq i\leq r/2$ and
$c_{r,r}=1$ independent
of $\alpha$ and $a$ such that 
\[
Z_{\alpha}(a[-1]^r\unit,q)=\left(\sum_{0\leq i\leq r/2}c_{r,r-2i}(a,\alpha)^{r-2i}
E_2(q)^{i}\right)q^{(\alpha,\alpha)/2}/\eta(q)^d.
\]
\end{lemma}

\begin{proof} We prove this by induction on $r.$ If $r=0$ the assertion
is clear. If $r=1$ the assertion follows from the fact
that $o(a[-1]\unit)=(a,\alpha)$ on $M(1)\otimes e^{\alpha}.$ Now we assume that
$r\geq 2.$ Using Corollary \ref{corollary:2.2.2} gives
\[
Z_{\alpha}(a[-1]^r\unit,q)=(a,\alpha)Z_{\alpha}(a[-1]^{r-1}\unit,q)
+(r-1){E}_2(q)Z_{\alpha}(a[-1]^{r-2}\unit,q).
\]
By induction, both $Z_{\alpha}(a[-1]^{r-1}\unit,q)$
and $Z_{\alpha}(a[-1]^{r-2}\unit,q)$ converge to  holomorphic functions
in $H,$ whence so does $Z_{\alpha}(a[-1]^r\unit,q).$ Indeed,
\begin{align*}
Z_{\alpha}(a[-1]^r\unit,q)&=(a,\alpha)\left(\sum_{0\leq i\leq (r-1)/2}c_{r-1,r-1-2i}(a,\alpha)^{r-1-2i}
E_2(q)^{i}\right)q^{(\alpha,\alpha)/2}/\eta(q)^d\\
&+(r-1)E_2(q)\left(\sum_{0\leq i\leq (r-2)/2}c_{r-2,r-2-2i}(a,\alpha)^{r-2-2i}
E_2(q)^{i}\right)q^{(\alpha,\alpha)/2}/\eta(q)^d\\
&=\left(\sum_{0\leq i\leq r/2}c_{r,r-2i}(a,\alpha)^{r-2i}
E_2(q)^{i}\right)q^{(\alpha,\alpha)/2}/\eta(q)^d.
\end{align*}
It is easy to express $c_{r,r-2i}$ in terms of the $c_{r-1,r-2j}$ and 
$c_{r-2,r-2k}$, and so obtain a recursive formula. We leave this 
to the reader. 
Since $c_{r-1,r-2j}$ and 
$c_{r-2,r-2k}$ are independent of $a,\alpha,$  $c_{r,r-2i}$ is also
independent of $a,\alpha.$
\end{proof}

For convenience, we set 
\[
f_{a,\alpha,r}(q)=\sum_{0\leq i\leq r/2}c_{r,r-2i}(a,\alpha)^{r-2i}
E_2(q)^{i}.
\]

Recall that $\{h_1,...,h_d\}$ is an orthonormal basis of $\frakh.$ Let $r_1,...,r_d$ be non-negative integers. Then by Corollary \ref{corollary:2.2.2}
and Lemma \ref{l3.2} we obtain
\begin{lemma} Let $v=h_1[-1]^{r_1}\cdots h_d[-1]^{r_d}\unit.$ Then 
$Z_{\alpha}(v,q)$ converges to a holomorphic function in the upper half plane
and
\[
Z_{\alpha}(v,q)=f_{h_1,\alpha,r_1}(q)\cdots f_{h_d,\alpha,r_d}(q)q^{(\alpha,\alpha)/2}/\eta(q)^d.
\]
\end{lemma}

Recall that $V_{L+\l_i}=\oplus_{\alpha\in L+\l_i}M(1)\otimes e^{\alpha}.$ 
The following corollary is immediate:
\begin{corollary}\label{c3.3} Let $v$ be as before. 
The function $Z_i(v,q)$ is equal to
\[
(\sum_{\alpha\in L+\l_i}f_{h_1,\alpha,r_1}(q)\cdots f_{h_d,\alpha,r_d}(q)q^{(\alpha,\alpha)/2})/\eta(q)^d.
\]
\end{corollary}

\begin{theorem}\label{t3.4} Let $v\in M(1).$ Then $Z_i(v,q)$ can be 
expressed as a finite sum
\[
Z_i(v,q)=\sum_{j}f_j(q)\Theta_{L+\l_i,k_j}(a_j,q)/\eta(q)^d
\]
where $f_j(q)\in Q,$ $k_j\geq 0,$ $a_j\in\frakh$ and
\[
\Theta_{L+\l_i,k_j}(a_j,q)=\sum_{\alpha\in L+\l_i}(a_j,\alpha)^{k_j}q^{(\alpha,\alpha)/2}.
\]
\end{theorem}

\begin{proof} First we take $v=a[-n]^{r_n}\cdots a[-1]^{r_1}\unit$
for $a\in\frakh$ with $(a,a)=1.$ Then from  Corollary \ref{corollary:2.2.2}
we have
\[
Z_i(v,q)=\sum_{s\geq 0}c_sg_s(q)Z_i(a[-1]^s\unit,q)
\]
for some scalars $c_s$ and quasi-modular forms $g_s(q)\in Q.$ Now the 
result follows from Corollary \ref{c3.3}. For arbitrary $v$ we can assume that
$v$ is a monomial in $h_j(-n)$ for $j=1,\cdots,d$ and $n>0.$ 
Use the result for 
$a[-n]^{r_n}\cdots a[-1]^{r_1}\unit$ and Corollary \ref{corollary:2.2.2} again
to see  that $Z_i(v,q)$
is a linear combination of functions of the form
\[
f(q)\Theta_{L+\l_i,t_1,...,t_d}(h_1,...,h_d,q)/\eta(q)^d
\]
where $f(q)$ belong to $Q,$ $t_j$ are non-negative integers 
and 
\[
\Theta_{L+\l_i,t_1,...,t_d}(h_1,...,h_d,q)=\sum_{\alpha\in L+\l_i}(h_1,\alpha)^{t_1}\cdots
(h_d,\alpha)^{t_d}q^{(\alpha,\alpha)/2}.
\]
It is easy to see that $\Theta_{L\!+\!\l_i,t_1,...t_d}(h_1,...,h_d,q)$
is a linear combination of $\Theta_{L+\l_i,t_1\!+\!\cdots\!+\!t_d}(b,q)$
for $b\in\frakh.$ This completes the proof of the theorem.
\end{proof}

    We now concentrate on the case $i = 0.$ Note that $Z_0(v, q) = Z(v, 
q).$ We prove
\begin{theorem}\label{t4.2.6}
Let $Q(N, \epsilon)$ be as in (\ref{1.0}) in the introduction. 
For each $v$ in $V_L,$ the  function $Z(v, q)\eta(q)^d$ lies in $Q(N, 
\epsilon).$
\end{theorem}
\begin{proof} After theorem \ref{t3.4} it suffices to show that if $a$ in 
$\frakh$ satisfies $(a, a) = 1$ then the function $\Theta_{L,k}(a, q)$ lies 
in $Q(N, \epsilon).$ This is essentially established in the paper \cite{DM2}. 
In the notation of \cite{DM2}, the present function $\Theta_{L, k}(a, q)$ is 
denoted  $\theta(Q, a, k, \tau)$ where $Q$ is the integral quadratic form 
which corresponds to $L.$ If $k$ is odd then the function is identically 
zero. If $k = 2l$ is even then Theorem 2 of (loc cit) says that the 
function
\begin{equation}\label{4.2.1}
\psi(Q, a, 2l, \tau)  =  \sum_{t=0}^l\gamma(t,2l)E_2(\tau)^t 
\theta(Q, a, 2l - 2t, \tau)
\end{equation}
is a holomorphic modular form in $M(N, \epsilon).$ In (\ref{4.2.1}), 
$\gamma(t, 2l)$ is a certain non-zero constant equal to $1$ when $t = 0.$ Moreover $\theta(Q, a, 0, \tau)$ is just the theta function of $L.$ Using this 
information, it follows from (\ref{4.2.1}) and induction on $l$ that 
$\theta(Q,a, 2l, \tau)$ indeed lies in $M(N, \epsilon),$ as required.
\end{proof}

    We are now in a position to complete the proof of Theorem \ref{thmm2}. In 
Theorem \ref{t4.2.6} we have shown that each trace function $Z(v, q)$ 
satisfies (\ref{1.3}) with $f(v, q)$ some element of $Q(N, \epsilon).$ We have 
to show that in fact $f(v, q)$ lies in $M(N, \epsilon).$ First note that 
we may assume without loss that $v$ lies in $(V_L)_{[k]}$  for some integer $k$ 
(cf. section 2.1). This puts us in a position to apply Proposition
\ref{pa}, which tell us that $Z(v, q)$ is a generalized modular form in 
the sense of the introduction. Thus the following result completes 
the proof of Theorem \ref{thmm2}:
\begin{proposition}\label{p4.2.7} An element of $Q(N, \epsilon)$ is a modular 
form of weight $k$ if, and only if, it is a generalized modular form of weight 
$k.$
\end{proposition}
\begin{proof} It is enough to prove sufficiency. Let $f(\tau)\in Q(N, \epsilon),$ of weight $k,$ be expressed in the form
\begin{equation}\label{4.2.2}
                f(\tau)  =\sum_{i=0}^m f_i (\tau) E_2(\tau)^i
\end{equation}
with each  $f_i (\tau)$ a form in $M(N, \epsilon)$ of weight $k - 2i.$ As 
$f(\tau)$ is also a generalized modular form of weight $k,$ each of its 
$\gamma$-transforms ($\gamma\in SL(2,\Z)$) has a $q$-expansion. In particular,
the $S$-transform of $f(\tau)$  yields an equality of the shape
\begin{equation}\label{4.2.3}
   \tau^{- k} f(S \tau)  =\sum_{n}b(n) q^{n/N}.
\end{equation}
   The $\gamma$-transform of $E_2(\tau)$ is well-known. In particular, we have
\begin{equation}\label{4.2.4}
 E_2(S \tau)  =  \tau^2 E_2(\tau) - \tau/2 \pi i.
\end{equation}
Comparing (\ref{4.2.2}) - (\ref{4.2.4}) we obtain an equality
\begin{equation}\label{4.2.5} 
\sum_n b(n) q^{n/N}  =\sum_{i=0}^m \tau^{- k + 2i} f_i (S 
\tau) ( E_2(\tau) - 1/\log q)^i.
\end{equation}
Since each $\tau^{-k + 2i}f_i (S \tau)$ has a $q$-expansion, an equality 
like (\ref{4.2.5}) can only hold if all $f_i (\tau)$  are zero for $i>0.$ 
Thus  $f(\tau)  =  f_0(\tau)$  is indeed modular.
\end{proof}

    This completes the proof of Proposition \ref{p4.2.7}, and hence that of 
Theorem \ref{thmm2}, with the following caveat: in Theorem \ref{t4.2.6} and 
Proposition \ref{p4.2.7} we were implicitly dealing with  modular forms of 
integral  weight. But it is easy to see that these results and those 
of \cite{DM2} extend to the half-integral case i.e., the case in which 
the lattice $L$ has odd rank. We leave the details to the interested 
reader.

\section{Spherical harmonics and trace functions}

In this section we first discuss the relation between spherical 
harmonics and highest weight vectors for the Virasoro algebra.
Each  spherical 
harmonic $P$ of degree $k$  gives rise to a highest weight vector
$v_P$ of weight $k$ in an canonical way. We then prove 
Theorem \ref{thmm3}.

\subsection{Spherical harmonics and primary vectors}

We continue our discussion of the vertex operator algebra $V_L.$ In particular,
$M(1)$ is a vertex operator subalgebra of $V_L.$ 

An element in $V_L$ is called a \textit{primary state} 
(or singular vector or highest
weight vector) if it satisfies $L(n)v = 0$ for all
$n\in\Z_+$. Thanks to the Virasoro commutator relations, this is
equivalent to $L(1)v = L(2)v = 0$.

\begin{lemma} Let $v\in M(1)$ be a polynomial in the variables
$h_i(-1), \,(1\leq i\leq d)$. Then $v$ is quasiprimary, i.e., $L(1)v =
0$.
\end{lemma}
\begin{proof} Recall that $\{h_1,...,h_d\}$ is an orthonormal basis of
$\frakh.$ Note that
\[ L(1)=\frac{1}{2}\sum_{i=1}^d\sum_{k\in\Z}\NO h_i(1-k)h_i(k)\NO
=\sum_{i=1}^d\sum_{k\geq 1}h_i(1-k)h_i(k)
\] and both $h_i(k)$ for $k\geq 2$ and $h_i(0)$ annihilate $v$. Since all
summands of $L(1)$ contain one of these as a factor, we immediately see
that $L(1)v = 0$.
\end{proof}

\begin{lemma} Let $v\in M(1)$ be a polynomial in the variables
$h_i(-1), \,(1\leq i\leq d)$. Then $v$ is primary if, and only if,
$v$ is a spherical harmonic, i.e., $\Delta v = 0$ where
\[
\Delta =
\sum_{i=1}^d \frac{\partial^2}{\partial h_i(-1)^2}.
\]
\end{lemma}
\begin{proof} It suffices to consider the condition $L(2)v = 0$. Since
\[ L(2) = \frac{1}{2}\sum_{i=1}^d h_i(1)^2 + \sum_{i=1}^d
\sum_{k\geq 2}h_i(2-k)h_i(k)
\] and $h_i(k)$ annihilates $v$ if $k\geq 2,$ we see that $L(2)v= 0$ is
equivalent to $\sum_{i=1}^d h_i(1)^2v = 0$. Now since
$h_i(1)h_i(-1)^n\unit = nh_i(-1)^{n-1}\unit$, $h_i(1)$ acts as
$\partial/\partial h_i(-1)$ on $M(1)$.
\end{proof}

\subsection{Proof of Theorem 3}

Let $P$ be a spherical harmonic of degree $k$:
\[
P=P(x_1,\dots,x_d),\ \Delta P= 0
\]
where $\Delta =\sum_{i=1}^d \partial^2/\partial x_i^2$. 
Let $v_P$ be the corresponding primary state
in $M(1)$, i.e.,
\[
v_P= P(h_1(-1),\ldots,h_d(-1))\unit.
\]

We note that $o(v_P) = v_P(k-1).$ Now let $P =
\sum_{I}c_Ix^I$ where $c_I\in\C$ and $x^I = x_1^{a_1}\cdots
x_d^{a_d}$ and where $I = (a_1,\ldots,a_d)$ is some ordered tuple
of non-negative integers such that $\sum a_i = k$. Then from Section 3.1,
\[ Y(v_P,z) = \sum_{I}c_IY(v_I,z) = \sum_{I}c_I\NO h_I(z)\NO
\] where $h_I(z) = \prod_{i=1}^d h_i(z)^{a_i}$ and $h_i(z) =
\sum_{n\in\Z}h_i(n)z^{-n-1}$. Therefore
\begin{equation}\label{5.2.1}
o(v_P) =\sum_{I}c_I\sum\NO h_1(n_{11})\cdots
h_1(n_{1a_1})h_2(n_{21})\cdots h_2(n_{2a_2})\cdots h_d(n_{d
1})\cdots h_1(n_{d a_d})\NO
\end{equation}
where the indices in the second sum range over all integers such that
$\sum_{i=1}^d
\sum_{j=1}^{a_i}n_{ij} = 0$. In particular, this includes  the case with
all
$n_{ij} = 0$:
\[
\sigma = \sum_{J}c_Jh_1(0)^{a_1}\cdots h_d(0)^{a_d} =
P(h_1(0),\dots, h_d(0)).
\] As $h_i(0)$ operates on $u\otimes e^\alpha$ as in (\ref{4.1.1}),
we see that
\[
\sigma:u\otimes e^\alpha\longmapsto \left(\sum_{J}c_J (h_1,\alpha)^{a_1}
\cdots (h_d,\alpha)^{a_d}\right) u\otimes e^\alpha.
\] But this is exactly $P(\alpha)u\otimes e^\alpha$. It follows that
\[
\tr_{V_L}\sigma q^{L_0-\frac{d}{24}} = \sum_{\alpha\in L}P(\alpha)
q^{\frac{(\alpha,\alpha)}{2}}/\eta(\tau)^d =
\frac{\theta_L(P,\tau)}{\eta(\tau)^d}.
\] To prove the theorem, it remains to show that all other terms in (\ref{5.2.1}) 
have trace $0$ in their action on $V_L.$  To see this, we adopt a different approach and note that to prove the theorem it is enough to establish
it for a
spanning set of spherical harmonics of degree $k$. These are given, for
example [H], by
\[ P = (t_1x_1+\cdots+t_d x_d)^k
\] for $(t_1,\dots,t_d)\in\C^d$ such that $\sum_{i=1}^d t_i^2 =
0$. Set $h_0 = \sum_{i=1}^d t_ih_i\in\frakh$, and note that $(h_0,h_0) =
\sum_{i=1}^d t_i^2 = 0$. Then
\[ Y(v_P,z) = \NO (\sum_{i=1}^d t_i h_i(z))^k\NO = \NO h_0(z)^k\NO.
\] Clearly
\[ o(v_P) = h_0(0)^k + \sum_{ n\neq
0}\binom{k}{2}h_0(0)^{k-2}h_0(n)h_0(-n)+\cdots.
\] Then it is enough to show that each operator of the form
\begin{equation}\label{5.2.2}
U= h_0(n_1)h_0(n_2)\cdots h_0(n_r),\quad (r\geq 2,\, not\ all\ n_i \ are\ 
zero,\ \sum n_i = 0)
\end{equation}
 has trace $0$ on $M(1)$. Let us take a basis
$h_0,k_0,k_1,\ldots,k_{d-2}$ of $\frakh$ with
\[ (h_0,k_0) = 1, \quad (h_0,k_i) = 0,\quad (1\leq
i\leq d-2).
\] Then $h_0(m)$ commutes with all $h_0(m), k_i(m),\,(1\leq i\leq
d,m\in\Z)$. Now consider a vector in $M(1)$ of the form $p\cdot v $
where $p$ is a polynomial in $k_0(-m),m\in\Z_+$ and $v$ does not
involve any mode $k_0(m), m\in \Z$. Since some $n_i$ occurring in
(\ref{5.2.2}) is positive, we see that the operator $U$ maps $pv$ to
$0$ or $p' v'$ where $p'$ has lower degree than $p$.  Therefore
the trace of $U$ on $M(1)$ is $0$ as the operator $U$ never has
non-zero eigenvalues. This completes the proof of Theorem \ref{thmm3}.

As explained in the introduction, the result of Waldspurger [W] together
with Theorem \ref{thmm3} implies the
\begin{corollary} Let $f(q)$ be any cusp-form of level $1.$ There there is a primary field $v$ in the lattice VOA $V_{E_8}$ such that $Z(v,q)=f(q)/\eta(q)^8.$
\end{corollary}

\end{document}